\documentclass{amsart}
\usepackage{amsmath}
\usepackage{amssymb}
\usepackage{amsthm}
\begin{document}
\def\eq#1{{\rm(\ref{#1})}}
\theoremstyle{plain}
\newtheorem{thm}{Theorem}[section]
\newtheorem{lem}[thm]{Lemma}
\newtheorem{prop}[thm]{Proposition}
\newtheorem{cor}[thm]{Corollary}
\theoremstyle{definition}
\newtheorem{dfn}[thm]{Definition}
\newtheorem{rem}[thm]{Remark}
\def\Ker{\mathop{\rm Ker}}
\def\Coker{\mathop{\rm Coker}}
\def\ind{\mathop{\rm ind}}
\def\Re{\mathop{\rm Re}}
\def\vol{\mathop{\rm vol}}
\def\SO{\mathbin{\rm SO}}
\def\Im{\mathop{\rm Im}}
\def\min{\mathop{\rm min}}
\def\Spec{\mathop{\rm Spec}\nolimits}
\def\Hol{{\textstyle\mathop{\rm Hol}}}
\def\ge{\geqslant}
\def\le{\leqslant}
\def\Z{{\mathbin{\mathbb Z}}}
\def\R{{\mathbin{\mathbb R}}}
\def\N{{\mathbin{\mathbb N}}}
\def\al{\alpha}
\def\be{\beta}
\def\ga{\gamma}
\def\de{\delta}
\def\ep{\epsilon}
\def\io{\iota}
\def\ka{\kappa}
\def\la{\lambda}
\def\ze{\zeta}
\def\th{\theta}
\def\vp{\varphi}
\def\si{\sigma}
\def\up{\upsilon}
\def\om{\omega}
\def\De{\Delta}
\def\Ga{\Gamma}
\def\Th{\Theta}
\def\La{\Lambda}
\def\Om{\Omega}
\def\ts{\textstyle}
\def\sst{\scriptscriptstyle}
\def\sm{\setminus}
\def\op{\oplus}
\def\ot{\otimes}
\def\bigop{\bigoplus}
\def\iy{\infty}
\def\ra{\rightarrow}
\def\longra{\longrightarrow}
\def\dashra{\dashrightarrow}
\def\t{\times}
\def\w{\wedge}
\def\d{{\rm d}}
\def\bs{\boldsymbol}
\def\ci{\circ}
\def\ti{\tilde}
\def\ov{\overline}
\def\md#1{\vert #1 \vert}
\def\nm#1{\Vert #1 \Vert}
\def\bmd#1{\big\vert #1 \big\vert}
\def\cnm#1#2{\Vert #1 \Vert_{C^{#2}}} 
\def\lnm#1#2{\Vert #1 \Vert_{L^{#2}}} 
\def\bnm#1{\bigl\Vert #1 \bigr\Vert}
\def\bcnm#1#2{\bigl\Vert #1 \bigr\Vert_{C^{#2}}} 
\def\blnm#1#2{\bigl\Vert #1 \bigr\Vert_{L^{#2}}} 
\title[Deformations of SLag Submanifolds; Fredholm Alternative Approach]{Deformations of Special Lagrangian Submanifolds; An Approach via Fredholm Alternative}
\author{Sema Salur}
\address {Department of Mathematics, Northwestern University, IL, 60208 }
\email{salur@math.northwestern.edu }

\begin{abstract} In an earlier paper, \cite{salur}, we showed that
the moduli space of deformations of a smooth, compact, orientable
special Lagrangian submanifold $L$ in a symplectic manifold $X$
with a non-integrable almost complex structure is a smooth
manifold of dimension ${H}^1(L)$, the space of harmonic 1-forms on
$L$. We proved this first by showing that the linearized operator
for the deformation map is surjective and then applying the Banach
space implicit function theorem. In this paper, we obtain the same
surjectivity result by using a different method, the Fredholm
Alternative, which is a powerful tool for compact operators in
linear functional analysis.
\end{abstract}
\date{}
\maketitle
\section{Introduction}
\label{co1}

In \cite{mclean}, McLean showed that the moduli space of nearby
submanifold of a smooth, compact special Lagrangian submanifold
$L$ in a Calabi-Yau manifold $X$ is a smooth manifold and its
dimension is equal to the dimension of ${H}^1(L)$, the space of
harmonic 1-forms on $L$. Special Lagrangian submanifolds have
attracted much attention after Strominger, Yau and Zaslow proposed
a mirror Calabi-Yau construction using special Lagrangian
fibration \cite{syz}. For more information about special Lagrangian submanifolds and examples, see \cite{bryant1}, \cite{bryant2}, \cite{HL}.

One can also define special Lagrangian submanifolds of symplectic
manifolds equipped with a nowhere vanishing complex valued
$(n,0)$-form, \cite{salur}. Such symplectic manifolds were studied
recently by Smith, Thomas and Yau in \cite{sty}.

In \cite{salur}, we showed that the moduli space of special
Lagrangian deformations of $L$ in a symplectic manifold with
non-integrable almost complex structure is also a smooth manifold
of dimension $b_1(L)$, the first Betti number of $L$. In order to
prove this result we first modified the definition of a special
Lagrangian submanifold for symplectic manifolds, extended the
parameter space of deformations and showed that the linearization
of the deformation map is onto and finally applied the infinite
dimensional Banach space implicit function theorem.


In this paper, we obtain the same result by a different approach. In particular, we show that the linearized operator for the deformation map is invertible by
using Fredholm Alternative, a technique from linear functional analysis.

\section{Deformations of Special Lagrangian Submanifolds}

Let $(M^{2n},\om, J, g,\Om )$ be a Calabi-Yau manifold with a
K\"{a}hler 2-form $\om$, a complex structure $J$, a compatible
Riemannian metric $g$ and a nowhere vanishing holomorphic
$(n,0)$-form $\Om $. Then one can define a special Lagrangian
submanifold of $M$.

\begin{dfn} An $n$-dimensional submanifold $L\subseteq M$ is
{\em special Lagrangian} if $L$ is Lagrangian (i.e. $\om|_L\equiv
0$) and $Im(\Om)$ restricts to zero on $L$. Equivalently,
Re$(\Om)$ restricts to be the volume form on $L$ with respect to
the induced metric.
\end{dfn}

McLean studied the deformations of compact special Lagrangian
submanifolds in Calabi-Yau manifolds and proved the following
theorem, \cite{mclean}.

\begin{thm} The moduli space of all deformations of a smooth, compact,
orientable special Lagrangian submanifold $L$ in a Calabi-Yau
manifold $M$ within the class of special Lagrangian submanifolds
is a smooth manifold of dimension equal to dim$(H^1(L))$.
\label{co2thm}
\end{thm}


One natural generalization of McLean's result is for symplectic
manifolds. Now let $(X,\omega,J,g,\xi)$ denote a $2n$-dimensional
symplectic manifold $X$ with symplectic 2-form $\omega$, an almost
complex structure $J$ which is tamed by $\omega$, the compatible
Riemannian metric $g$ and a nowhere vanishing complex valued
$(n,0)$-form $\xi=\mu+i\beta$, where $\mu$ and $\beta$ are real
valued $n$-forms.




Note that the holomorphic form $\Omega$ is a closed form on
Calabi-Yau manifolds. However, on a symplectic manifold which is
equipped with a nowhere vanishing complex-valued $(n,0)$ form
$\xi$, it is not necessarily closed.



For more general special Lagrangian calibrations, one can
introduce an additional term $e^{i\theta}$, where for each fixed
angle $\theta$ we have a corresponding form $e^{i\theta}\xi$ and
its associated geometry. Here $\theta$ is the phase factor of the
calibration and using this as the new parameter one can define
special Lagrangian submanifolds in a symplectic manifold and study
their deformations, \cite{salur}.




\begin{dfn} An $n$-dimensional submanifold $L\subseteq X$ is
{\em special Lagrangian} if $L$ is Lagrangian (i.e.
$\omega|_L\equiv 0$) and $Im(e^{i\theta}\xi)$ restricts to zero on
$L$, for some $\theta\in {\mathbb R}$. Equivalently,
$Re(e^{i\theta}\xi)$ restricts to be the volume form on $L$ with
respect to the induced metric.
\end{dfn}

Now we recall the basics of a technique for compact operators in
linear functional analysis, known as Fredholm Alternative. One can
find more information about the subject in \cite{evans} and
\cite{lax}.

Let $\mathcal X$ and $\mathcal Y$ be two real Banach spaces.

\begin{dfn}
A bounded linear operator  $\mathcal K$ : $\mathcal X \rightarrow
\mathcal Y$ is called compact provided for each bounded sequence
$\{u_k\}_{k=1}^{\infty}$ is precompact in $\mathcal Y$.
\end{dfn}

Now let $H$ denote a real Hilbert space, with inner product $<,>$.

\begin{thm} Let $\mathcal K:H\rightarrow
H$ be a compact operator. Then

i) $\ker(I-\mathcal K)$ is finite dimensional,

ii) Range$(I-\mathcal K)$ is closed,

iii) Range$(I-\mathcal K)=\ker(I-\mathcal K^*)^\bot$,

iv) $\ker(I-\mathcal K)=\{0\}$ if and only if Range$(I-\mathcal
K)=H$

v) dim $\ker(I-\mathcal K)=dim\ker(I-\mathcal K^*)$.
\label{specthm1}
\end{thm}

\begin{rem} Theorem \ref{specthm1} asserts in particular either

(a) for each $f\in H$, the equation $u-\mathcal Ku=f$ has a unique solution

\noindent or else

(b) the homogeneous equation $u-\mathcal Ku=0$ has solutions $u\neq 0$.

In addition, should (a) obtain the space of solutions of the
homogeneous problem is finite dimensional and the nonhomogeneous
equation $u-\mathcal Ku=f$ has a solution if and only if
$f\in\ker(I-\mathcal K^*)^\bot$.

\end{rem}

Now we prove the following theorem, \cite{salur}, using the
Fredholm Alternative:

\begin{thm} Let $L$ be a smooth, compact, orientable special Lagrangian submanifold of a symplectic manifold $X$.
Then the moduli space of all deformations of $L$ in $X$ within the
class of special Lagrangian submanifolds is a smooth manifold of
dimension ${H}^1(L)$.
\end{thm}


\begin{proof}
Given a domain $\Omega$, let $C^{k,\alpha}(\Omega)$ denote the
H\"{o}lder norms defined as

\begin{equation*}
C^{k,\alpha}(\Omega)=\{ f\in C^k(\Omega) |\;\;
[D^{\gamma}f]_{\alpha,\Omega}<\infty , |\gamma|\leqq k \}
\end{equation*}

\noindent where

\begin{equation*}
[f]_{\alpha,\Omega}= \mathop{Sup}\limits_{x,y\in \Omega,\; x\neq y
} \frac{dist(f(x),f(y))}{(dist(x,y))^\alpha}\;\;\;  {\text in}\;\;
\Omega.
\end{equation*}

Then for a small vector field $V$ and a scalar $\theta\in
\mathbb{R}$, we define the deformation map as follows,
\begin{equation*}
F: C^{1,\alpha}(\Gamma(N(L)))\times\mathbb{R}\rightarrow
C^{0,\alpha}(\Omega^2(L))\oplus C^{0,\alpha}(\Omega^n(L))
\end{equation*}

\begin{equation*}
F(V,\theta)=((\exp_V)^*(-\omega), (\exp_V)^*(Im(e^{i\theta}\xi)).
\end{equation*}


Here $N(L)$ denotes the normal bundle of $L$, $\Gamma(N(L))$ the
space of sections of the normal bundle, and $\Omega^2(L)$,
$\Omega^n(L)$ denote the differential $2$-forms and $n$-forms,
respectively.

Since the symplectic form $\omega$ is closed on $X$ and the
restriction of $Im(e^{i\theta}\xi)$ is a top dimensional form on
$L$ the image of the deformation map $F$ lies in the closed
$2$-forms and closed $n$-forms. So by Hodge decomposition we get

\begin{equation*}
F:C^{1,\alpha}(\Gamma(N(L)))\times\mathbb{R}\rightarrow
C^{0,\alpha}(d\Omega^1(L))\oplus
C^{0,\alpha}(d\Omega^{n-1}(L)\oplus\mathcal{H}^n(L))
\end{equation*}

\noindent where $d\Omega^{n-1}(L)$ denotes the space of exact
$n$-forms and $\mathcal{H}^n(L)$ denotes the space of harmonic
$n$-forms on $L$.



In \cite{salur}, we computed the linearization of $F$ at (0,0),

\begin{equation*}
dF(0,0):C^{1,\alpha}(\Gamma(N(L)))\times\mathbb{R}\rightarrow
C^{0,\alpha}(d\Omega^1(L))\oplus
C^{0,\alpha}(d\Omega^{n-1}(L)\oplus\mathcal{H}^n(L))
\end{equation*}

\noindent where

\begin{equation*}
dF(0,0)(V,\theta)=\frac{\displaystyle\partial}{\displaystyle\partial{t}}F(tV,s\theta)|_{t=0,
s=0}
+\frac{\displaystyle\partial}{\displaystyle\partial{s}}F(tV,s\theta)|_{t=0,
s=0}
\end{equation*}






\begin{equation*}
=(-d(i_V\omega)|_L,(i_Vd\beta +d(i_V\beta))|_L +\theta)
\end{equation*}

\begin{equation*}
=(dv, \zeta+d*v+\theta), \;\;\; \text{where}\;\;
\zeta=i_V(d\beta)|_L.
\end{equation*}

\noindent Here $i_V$ is the interior derivative and $v$ is the
dual $1$-form to the vector field $V$ with respect to the induced
metric. For the details of local calculations see \cite{mclean},
\cite{salur}.

Let $x_1, x_2,...,x_{n}$ and $x_1, x_2,...,x_{2n}$ be the local
coordinates on $L$ and $X$, respectively. Then for any given
normal vector field $V=(V_1\frac {\partial}{\partial
x_{n+1}},...,V_n\frac {\partial}{\partial x_{2n}})$ to $L$ we can
show that

\begin{equation*}
\zeta=i_V(d\beta)|_L=-n(V_1\cdot g_1+...+V_n\cdot g_n)\text{dvol}
\end{equation*}

\noindent where $g_i$ $(0<i\leq n)$ are combinations of
coefficient functions in the connection-one forms.


One can decompose the $n$-form $\zeta=da+d^*b+h_2$ by using
Hodge Theory and because $\zeta$ is a top dimensional form on $L$,
$\zeta$ is closed and the equation becomes
\begin{equation*}
dF(0,0)(V,\theta) =(dv,da+d*v+h_2+\theta)
\end{equation*}

\noindent for some $(n-1)$-form $a$ and harmonic $n$-form $h_2$.
Also the harmonic projection for
$\zeta=-n(V_1.g_1+...+V_n.g_n)$dvol is given by
$(\displaystyle\int_L{-n(V_1.g_1+...+V_n.g_n)}$dvol)dvol and
therefore one can show that

\begin{equation*}
da=-n(V_1.g_1+...+V_n.g_n)\text{dvol}+(n\int_L{(V_1.g_1+...+V_n.g_n)}\text{dvol)dvol}
\end{equation*}
\noindent and
\begin{equation*}
h_2=(-n\int_L{(V_1.g_1+...+V_n.g_n)}\text{dvol)dvol}.
\end{equation*}

One should note that the differential forms $a$ and $h_2$ both depend on $V$.






The Implicit Function Theorem says that $F^{-1}(0,0)$ is a
manifold and its tangent space at $(0,0)$ can be identified with
the kernel of $dF$.

\begin{equation*}
(dv)\oplus(\zeta+d*v+\theta)=(0,0)
\end{equation*}

\noindent which implies

\begin{equation*}
dv=0 \;\;\; \text{and} \;\;\;
\zeta+d*v+\theta=da+d*v+h_2+\theta=0.
\end{equation*}

The space of harmonic $n$-forms $\mathcal{H}^n(L)$, and the space
of exact $n$-forms $d\Omega^{n-1}(L)$, on $L$ are orthogonal
vector spaces by Hodge Theory. Therefore, $dv=0$ and
$da+d*v+h_2+\theta=0$ is equivalent to $dv=0$ and $d*v+da=0$ and
$h_2+\theta=0$.

One can see that the special Lagrangian deformations (the kernel
of $dF$) can be identified with the $1$-forms on $L$ which satisfy
the following equations:

\begin{equation*}
\begin{split}
(i)\;\;&dv=0\\
(ii)\;\;&d*(v+\kappa(v))=0\\
(iii)\;\;&h_2+\theta=0.
\end{split}
\end{equation*}

Here, $\kappa(v)$ is a linear functional that depends on $v$ and
$h_2$ is the harmonic part of $\zeta$ which also depends on $v$.
These equations can be formulated in a slightly different way in
terms of decompositions of $v$ and $*a$.

If $v=dp+d^*q+h_1$ and $*a=dm+d^*n+h_3$ then we have

\begin{equation*}
\begin{split}
(i)\;\;&dd^*q=0\\
(ii)\;\;&\Delta(p\pm m)=0\\
(iii)\;\;&h_2+\theta=0.
\end{split}
\end{equation*}

This formulation of the solutions will provide the proof of the
surjectivity of the linearized operator without using $\kappa(v)$.


Now we show that the linearized operator is surjective at $(0,0)$.
Recall that the deformation map is given as


\begin{equation*}
F: C^{1,\alpha}(\Gamma(N(L)))\times\mathbb{R}\rightarrow
C^{0,\alpha}(d\Omega^1(L))\oplus
C^{0,\alpha}(d\Omega^{n-1}(L)\oplus\mathcal{H}^n(L)).
\end{equation*}









Therefore, for any given exact $2$-form $x$ and closed $n$-form
$y=u+z$ in the image of the deformation map (here $u$ is the exact
part and $z$ is the harmonic part of $y$), we need to show that
there exists a 1-form $v$ and a constant $\theta$ that satisfy the
equations,

\begin{equation*}
\begin{split}
(i)\;\; &dv=x\\
(ii)\;\; &d*(v+\kappa (v))=u\\
(iii)\;\;  &h_2+\theta =z.
\end{split}
\end{equation*}

\vspace{.1in}

Alternatively, we can solve the following equations for $p,q$ and
$\theta$.

\vspace{.1in}
\begin{equation*}
\begin{split}
(i) \;\; &dd^*q=x\\
(ii)\;\; & \Delta(p\pm m)=*u\\
(iii) \;\; & h_2+\theta=z.
\end{split}
\end{equation*}

\noindent where the star operator $*$ in $(ii)$ is defined on $L$.

For (i), since $x$ is an exact $2$-form we can write
$x=d(dr+d^*s+$harmonic form) by Hodge Theory. Then one can solve
(i) for $q$ by setting $q=s$.

For (ii), since $\Delta m=d^*dm=d^**a=*d**a=\pm *da$,
\begin{equation*}
\Delta(p\pm m)= \Delta p \pm *da
\end{equation*}

\noindent where $a$ depends on $p$ and we obtain

\vspace{.1in}
\begin{equation}
\Delta p
\pm(-n(V_1.g_1+...+V_n.g_n)+(n\displaystyle\int_L{(V_1.g_1+...+V_n.g_n)}\text{dvol}
))=*u
\label{eq1}
\end{equation}

We can show the solvability as follows: Since $V=(V_1,...,V_n)$ is
the dual vector field of the one form $v=dp+d^*q+h_1$ we can write
the equation \eq{eq1} as

\begin{equation}
\Delta p \pm (-n(v\cdot g)+(n{\displaystyle\int_L}{(v\cdot
g)}\text{dvol}))=*u
\end{equation}

\begin{equation}
\Delta p \pm (-n(dp+d^*q+h_1)\cdot
g+(n{\displaystyle\int_L}{(dp+d^*q+h_1)\cdot g}\text{dvol}))=*u.
\end{equation}

\noindent where $v\cdot g $ represents the action of the one form $v$ on the vector field $g=(g_1,..,g_n)$ and $n{\displaystyle\int_L}{(dp+d^*q+h_1\cdot g)}$dvol  is the harmonic projection of $-n(dp+d^*q+h_1)\cdot g$.

Then we get
\begin{equation*}
\Delta p \pm n( -(dp\cdot g)+{\displaystyle\int_L}{dp \cdot
g}\;\text{dvol} ) =*u\mp n[-(d^*q+h_1)\cdot
g+{\displaystyle\int_L}{(d^*q+h_1)\cdot g}\;\text{dvol}].
\end{equation*}

\noindent For simplicity we put
\begin{equation*}
*u\mp n[-(d^*q+h_1)\cdot g+{\displaystyle\int_L}{(d^*q+h_1)\cdot
g}\;\text{dvol} ]=h.
\end{equation*}

Since ${\displaystyle\int_L} *u=0$ and
${\displaystyle\int_L}{(d^*q+h_1)\cdot g}$ dvol is equal to the
harmonic projection of $(d^*q+h_1)\cdot g$, we get
${\displaystyle\int_L} h=0 .$

Since $L$ is a compact manifold without boundary, by Stoke's
Theorem,
\begin{equation*}
{\displaystyle\int_L}{dp\cdot g}\;\text{dvol}=
-{\displaystyle\int_L}{p\cdot \text{div} g}\;\text{dvol}
\end{equation*}

\noindent and the equation becomes
\begin{equation*}
\Delta p \pm n( -(dp\cdot g)-{\displaystyle\int_L}{p\cdot
\text{div} g}\;\text{dvol} )=h.
\end{equation*}

Then by adding and subtracting $p$ from the equation

\begin{equation*}
(\Delta - Id)p = [\pm n( -(dp\cdot g)-{\displaystyle\int_L}{p\cdot
\text{div} g}\;\text{dvol}) - p + h]
\end{equation*}

\noindent and

\begin{equation*}
p=(\Delta - Id)^{-1} [.....]p +\overline{h}= {\mathcal
K}(p)+\overline{h}
\end{equation*}

\noindent where $\overline{h}=(\Delta - Id)^{-1}h$.

\noindent Since

\begin{equation*}
||(\Delta - Id)^{-1}{\displaystyle\int_L}{p\cdot \text{div}
g}||_{L_1^2}\leq C|{\displaystyle\int_L}{p\cdot \text{div} g}|\leq
C ||p||_{L^2}
\end{equation*}

\noindent ${\mathcal K}(p)$ is a compact operator which takes
bounded sets in ${L^2}$ to bounded sets in ${L_1^2}$. Also note
that we assumed here $1\notin $spec$(\Delta)$, and if this is not
the case then we can modify the above argument by adding and
subtracting $\lambda p$, $\lambda \notin $spec$(\Delta)$ from the
equation.

Next we show that the set of solutions of the equation

\begin{equation}
\Delta p \pm n( -(dp\cdot g)-{\displaystyle\int_L}{p\cdot
\text{div} g}\;\text{dvol} )=0
\end{equation}

\noindent is constant functions and therefore of dimension 1. Note
that this set of solutions also satisfy the equation
$(Id-{\mathcal K})(p)=0$.

Also note that ${\displaystyle\int_L}{p\cdot \text{div} g}$ dvol
is a constant which depends on $p$. We denote this as $C(p)$. At
maximum values of $p$, $\Delta p$ will be negative which imply that
$C(p)\leq 0$ and at minimum values of $p$, $\Delta p$ will be positive
which imply that $C(p)\geq 0$ so $C(p)$ should be zero. Then the
maximum principle holds for the equation $\Delta p \pm n(
-(dp\cdot g))=0$ and since $L$ is a compact manifold without
boundary the solutions of this equation are constant functions.
Hence the dimension of the kernel of $(Id-{\mathcal K})$ is one.

Next we find the kernel of $(Id-{\mathcal K^*})$.

\begin{equation}
{\displaystyle\int_L}(\Delta p \pm n( -(dp\cdot g)- \int_L p\cdot
\text{div} g))q(y)dy
\end{equation}

\begin{equation*}
={\displaystyle\int_L}p\Delta q(y) \pm
n{\displaystyle\int_L}-(dp\cdot g) q(y)dy-
n{\displaystyle\int_L}(\int_L p\cdot \text{div} g) q(y)dy.
\end{equation*}
\begin{equation*}
={\displaystyle\int_L}p(y)\Delta q \pm n{\displaystyle\int_L}+(p
\text{div} (g\cdot q)(y)dy- n{\displaystyle\int_L} p(x)\cdot
\text{div} g(x) \int_L q(y)dy dx
\end{equation*}
\begin{equation*}
={\displaystyle\int_L}p(y)\Delta q \pm n{\displaystyle\int_L}+(p
\text{div} (g\cdot q)(y)dy- n{\displaystyle\int_L} p(y)\cdot
\text{div} g(y) \int_L q(x)dx dy
\end{equation*}
\begin{equation*}
={\displaystyle\int_L}p(y)(\Delta q \pm n(+\text{div} (g\cdot
q)-\text{div} g \int_L q(x)dx) dy.
\end{equation*}

 Since we assumed that $1\notin $spec$(\Delta)$,
dim $\ker($Id$-{\mathcal K^*})(\Delta - $Id$)=$dim $ \ker
($Id$-{\mathcal K^*})$ and the kernel of $($Id$-{\mathcal K^*})$
is equivalent to the solution space of the equation

\begin{equation}
\Delta q \pm n(+\text{div} (g\cdot q)-\text{div} g \int_L
q(x)dx)=0.
\end{equation}

By Fredholm Alternative, the dimension of this kernel is 1 and one
can check that a constant function $q=1$ satisfies this equation,
therefore the kernel consists of constant functions. Moreover
these functions satisfy the compatibility condition
$\displaystyle\int h.q =0$.

Then by Fredholm Alternative, Theorem \ref{specthm1}, we can
conclude the existence of solutions of the equation

\begin{equation*}
\Delta p \pm
(-n(V_1.g_1+...+V_n.g_n)+(n{\int_L}{(V_1.g_1+...+V_n.g_n)}\;\text{dvol}))=*u
\end{equation*}

(iii) is straightforward.

It follows from \cite{salur} that the image of the deformation map
$F_1$ lies in $d\Omega^1(L)$ and the image of $F_2$ lies in
$d\Omega^{n-1}(L)\oplus\mathcal{H}^n(L)$. So we conclude that $dF$
is surjective at $(0,0)$. Also since both the index of $d+*d^*(v)$
and $d+*d^*(v+\kappa(v))$ are equal to $b_1(L)$ the dimension of
tangent space of special Lagrangian deformations in a symplectic
manifold is also $b_1(L)$, the first Betti number of $L$. By
infinite dimensional version of the implicit function theorem and
elliptic regularity, the moduli space of all deformations of $L$
within the class of special Lagrangian submanifolds is a smooth
manifold and has dimension $b_1(L)$.

\end{proof}

\begin{rem}
One can study the deformations of special Lagrangian submanifolds
in much more general settings. In a forthcoming paper we plan to
study these deformations using the techniques which we developed
recently for associative submanifolds of $G_2$ manifolds,
\cite{AkSa},\cite{AkSa1}.
\end{rem}

{\small{\it Acknowledgements.} The author is grateful to Peng Lu
for many useful conversations and valuable comments.
 }

\end{document}